\documentclass{amsart}
\usepackage[]{amsmath}
\usepackage{amsthm}
\usepackage{amsfonts}
\usepackage{amssymb}
\usepackage{tikz}
\usetikzlibrary{cd}

\theoremstyle{definition}
\newtheorem{theorem}{Theorem}
\newtheorem{prop}{Proposition}

\newtheorem{cor}{Corollary}

\newtheorem{defi}{Definition}
\newtheorem{example}{Example}

\newtheorem{obs}{Observation}

\newcommand{\C}{\mathbb{C}}
\newcommand{\Z}{\mathbb{Z}}

\renewcommand{\S}{\mathbb{S}}

\renewcommand{\sl}{\mathfrak{sl}}

\renewcommand{\d}{\mathfrak{d}}
\newcommand{\so}{\mathfrak{so}}

\newcommand{\g}{\mathfrak{g}}

\newcommand{\Spin}{\operatorname{Spin}}
\newcommand{\id}{\operatorname{Id}}

\title{Spin Representations and Binary Numbers}
\author{Henrik Winther}
\date{February 29, 2024}
\begin{document}
\begin{abstract}
	We consider a construction of the fundamental spin representations of the simple Lie algebras $\so(n)$ in terms of binary arithmetic of fixed width integers.
	This gives the spin matrices as a Lie subalgebra of a $\Z$--graded associative algebra (rather than the usual $\mathbb{N}$--filtered Clifford algebra).
	Our description gives a quick way to write down the spin matrices, and gives a way to encode some extra structure, such as the real structure which is invariant under the compact real form, for some $n$. 
	Additionally we can encode the spin representations combinatorially as (coloured) graphs.
\end{abstract}
\maketitle

\section{Introduction}
Finite dimensional representations of complex semisimple Lie algebras are characterized by their highest weight. 
Any irreducible representation can be constructed from tensor products of fundamental representations.
Thus it is important to have good ways to construct fundamental representations.
There are a few approaches to doing so, and Verma module quotients give a universal way, with the caveat of going through an infinite dimensional module that must be quotiented by an infinite dimensional submodule.
Thus one might wish to avoid this.

A common method is via the following observation:
\begin{obs}
	For the classical families $A_n,B_n,C_n,D_n$, we can produce all fundamental representations by taking exterior powers of a tautological representation $V$,
	\begin{equation}
		\Lambda V = \bigoplus_k \Lambda^k V = \bigoplus_\alpha V_\alpha 
	\end{equation}
	and picking out the highest component from each power.
	With one exception: The fundamental spin representations of $\so(n,\C)$.
\end{obs}
Thus we want to construct the spin representations. 
Of course, several constructions are known.
The standard way to do so is via Clifford algebras.
This can be turned into a method which however requires taking many iterated Kronecker products of matrices \cite[p. 11-12]{BFGK}.
In the present work we will provide another method to constructing the representation matrices of the fundamental spin representation, for which all matrix coefficients are precomputed.
\begin{theorem}
	Let $N>0$. The fundamental spin representation $(\S_{2N+1}, \so(2N+1,\C))$ is equivalent to $(\C \{1,0\}^N, \langle A^{\pm k}\circ A^{\mp (k-1)}| 0\le k<N \rangle)$ where $\C \{1,0\}^N$ is a free complex vector space generated by binary strings of fixed length $N$, (fixed width binary numbers), and $A^{\pm k}$ is the linear map
	\begin{equation}
		A^\pm_k: [m]_2 \mapsto \begin{cases}
			[m \pm 2^k]_2 \text{ if this operation would not require any carrying in binary arithmetic}\\
			0 \text{ else }
		\end{cases}
	\end{equation}
	and adding or subtracting fractions $2^{-1}$ is simply disregarded.
	By $\langle,\rangle$ we mean the Lie algebra generated by taking iterated Lie brackets.
	Moreover, the (decomposable) spin representation of $(\S_{2N}, \so(2N,\C))$ is equivalent to $(\C \{1,0\}^N, \langle A^{N-2}\circ(A^{N-1}+A^{-(N-1)} ) ,A^{\pm k}\circ A^{\mp (k-1)}| 0\le k<N-1 \rangle)$.
\end{theorem}
In this theorem we state generators, but see equations \eqref{eq:decompress1}, \eqref{eq:decompress2}, and \eqref{eq:decompress3}, as well as Theorem \ref{thm:evenspin}, for the complete list of operators without needing any extra brackets.
See also Section \ref{sec:compactspin} for the explicit operators of the compact real form.

\subsection{Acknowledgements}
This article is based upon work from COST Action CaLISTA CA21109 supported by COST (European Cooperation in Science and Technology).
This research was partially supported by the UiT Aurora project MASCOT.

\section{Arithmetic operators}
We will consider the set $S_N = \{ 0,1\}^N$ to be the set of binary numbers in $N$ digits.
Here we keep track also of trailing zeroes on the right, i.e. these are ``fixed width'' numbers.
For example,
\begin{equation}
	S_3 = \{ 000, 001, 010, 011, 100, 101, 110, 111\} = \{ [0]_2^3, [1]_2^3, [2^3,]_2^3, [3]_2^3, [4]_2^3, [5]_2^3, [6]_2^3, [7]_2^3, \}.
\end{equation}
We are going to consider functions on $S_N$ that are defined in terms of arithmetic formulas $m \mapsto m + k$ for $k\in \mathbb{Z}$.
However, we are going disregard carrying, and moreover, if it is not possible to perform the operation without carrying, then we define that our arithmetic function will simply fix its argument $m$.
\begin{example}
	We have the function $+2: m \mapsto m+2$, on $S_2$, given by
	\begin{equation}
		\begin{split}
			00 \mapsto 10\\
			01 \mapsto 11\\
			10 \mapsto 10\\
			11 \mapsto 11\\
		\end{split}
	\end{equation}
\end{example}
Of course, each integer $k$ defines a family of arithmetic functions in this way, one for each $S_N$.
Let us extend this to also include formal sums of arithmetic functions, such as $+k+l$.
We also define these to leave the argument unchanged whenever any carrying is necessary.
\begin{example}
	We have the function $+2-1: m \mapsto m+2-1$, on $S_2$, given by
	\begin{equation}
		\begin{split}
			00 \mapsto 00\\
			01 \mapsto 10\\
			10 \mapsto 10\\
			11 \mapsto 11\\
		\end{split}
	\end{equation}
\end{example}
Next we need to promote our arithmetic functions to linear operators.
To do this, we define 
\begin{equation}
	\S_{2N+1} = \C S_N = \C \{0,1\}^N
\end{equation}
as the free complex vector space over $S_N$.
We will eventually identify this with the a spin module, but in fact this vector space also has another name.
It is called the space of $n$-qubit states, in the context of quantum computing.
\begin{defi}
	For each arithmetic function $f: S_N\rightarrow S_N$, we define the linear operator
	\begin{equation}
		\bar f: \S_{2N+1}\rightarrow \S_{2N+1}	
	\end{equation}
	by
	\begin{equation}
		m \mapsto \begin{cases}
			0 \, \text{ if } f(m) = m\\
			f(m) \, \text{ else } 
		\end{cases}
	\end{equation}
	for $m\in S_N$, and extending by linearity to $\S_{2N+1}$.
	Operators defined in this way will be called \emph{arithmetic operators}.
\end{defi}
Each integer $k$ defines a family of arithmetic operators $\overline{+k}$, one for each $\S_{2N+1}$.

Let us introduce one more kind of operator.
\begin{defi}
	Let the $k$-parity operator be given by
	\begin{equation}
		m \mapsto p_k(m) = (-1)^{[m]_2(k)}m,
	\end{equation}
	where $[m]_2(k)$ is the $k$'th binary digit of the fixed-width number $m$, and extended to all $\S_{2N+1}$ by linearity.
\end{defi}
\begin{prop}\label{prop:parityprops}
	We have $p_k\circ p_k = \id$.
	The operator $p_k$ commutes with $\overline{\pm 2^l}$ for $l \not = k$.
	We have $p_k \circ \overline{\pm 2^k} = - \overline{\pm 2^k} \circ p_k$, and in particular,
	\begin{equation}\label{eq:pk-selfnormal}
		\begin{split}
			[p_k, \overline{2^k}] =  -2 \cdot \overline{ 2^k},\\
			[p_k, \overline{-2^k}] =  2 \cdot \overline{- 2^k},
		\end{split}
	\end{equation}
	and we have
	\begin{equation}\label{eq:pk-iscomm}
		[\overline{-2^k}, \overline{+2^k}] = p_k,
	\end{equation}
	and together \eqref{eq:pk-selfnormal} and \eqref{eq:pk-iscomm} shows that $(\overline{-2^k}, p_k, \overline{+2^k})$ is an $\mathfrak{sl}_2$-triple.
	We also have
	\begin{equation}
		[p_k,p_l] = 0
	\end{equation}
	for all $k,l$.
\end{prop}
\begin{proof}
	The arithmetic operator $\pm 2^{k}$ only changes the $k$'th bit of a given binary number. This does not change the $l$--parity unless $k=l$, and thus if $k$ and $l$ are different, they commute.
	In the case that $k=l$, it is easy to show that the commutation relations are as above.
\end{proof}
We can say more about relations between arithmetic operators.
\begin{prop}\label{prop:ariprops}
	We have 
	\begin{itemize}
		\item $[\overline{+2^n}, \overline{+2^m}] = 0$ for all $n,m$,
		\item $[\overline{-2^n}, \overline{-2^m}] = 0$ for all $n,m$,
		\item $[\overline{-2^n}, \overline{+2^m}] = 0$ for $m \not= n$.
	\end{itemize}
	If $0\le n<N$, then we have $\S_{2N+1} = \ker(\overline{2^n}) \oplus \ker(\overline{-2^n})$.
	Additionally, we have
	\begin{equation}\label{eq:compo}
		\begin{split}
			\overline{-2^{n-1}} \circ \overline{+2^{n-1}} = \tfrac{1}{2}(\id_{\S_{2N+1}} + p_n)\\
			\overline{+2^{n-1}} \circ \overline{-2^{n-1}} = \tfrac{1}{2}(\id_{\S_{2N+1}} - p_n)
		\end{split}
	\end{equation}
	Moreover,
	\begin{equation}
		\left[ [\, \overline{+2^{n-1}}, \overline{-2^{n-1}} \,] \right] = \id_{\S_{2N+1}},
	\end{equation}
	where $\left[ [,] \right]$ is the anti-commutator.
\end{prop}
\begin{proof}
	The first three statements follow because $(\overline{\pm 2^n})^2 = 0$, and if $m,n$ are different then the operators make non-trivial changes only on disjoint parts of the binary expansions of numbers.
	The equalities \eqref{eq:compo} can be obtained by comparing values:
	Suppose first that $[m]_2$ has a 1 in the $k$th position.
	Then $(\id_{\S_{2N+1}} + p_k)(m) = 0 =  \overline{-2^{n-1}} \circ \overline{+2^{n-1}}(m)$.
	But if $[m]_2$ has a 0 in the $k$th position, then
	$(\id_{\S_{2N+1}} + p_k)(m) = 2 =  2\cdot \overline{-2^{n-1}} \circ \overline{+2^{n-1}}(m)$.
	The next case is similar and the anticommutator follows.
\end{proof}

The next statement is an easy consequence of Prop. \ref{prop:parityprops} and Prop. \ref{prop:ariprops}:
\begin{prop}
	The Lie algebra generated by $\{\overline{-2^k}, p_{k+1}, \overline{+2^k} | k = 0 \dots N-1\}$ is isomorphic to
	\begin{equation}
		\mathfrak{sl}_2(\C)_1 \oplus \dots \oplus \mathfrak{sl}_2(\C)_N,
	\end{equation}
	and its representation on $\S_{2N+1}$ is equivalent to 
	\begin{equation}
		\C^2_1 \otimes \cdots \otimes \C^2_N,
	\end{equation}
	where the action of $\mathfrak{sl}_2(\C)_m$ is standard on $\C^2_m$ and trivial on $\C^2_l$ for $l\not=m$.
\end{prop}

\begin{prop}
	Let $k,l>0$ be natural numbers.
	Then if the binary expansions $[k]_2$ and $[l]_2$ have disjoint support,
	\begin{equation}
		\overline{\pm k} \circ \overline{\pm l} = \overline{\pm l} \circ \overline{\pm k} = \overline{\pm k \pm l},
	\end{equation}
	and
	\begin{equation}
		\overline{+k} \circ \overline{-l} = \overline{-l} \circ \overline{+k} = m \mapsto \begin{cases}
			m+k-l  \text{ if } [m]_2 \text{ has disjoint support from } [k+l]_2\\
			0 \text{ else}
		\end{cases}.
	\end{equation}
\end{prop}
\begin{proof}
	This clear from evaluating on arbitrary numbers.
\end{proof}

\begin{prop}\label{prop:zgrading}
	The algebra generated by arithmetic and parity operators comes equipped with a natural $\Z$-grading, by declaring $\overline{\pm k}$ to be an element of pure gradation $\pm k$, and parity operators $p_i$ to be of pure gradation zero.
\end{prop}
\begin{proof}
	This is well-defined because the algebra is generated by pure gradation elements, and a product of pure terms of gradations $\pm k_1,\dots \pm k_l $ will always map $m \in S_{N}$ to something proportional to $m \pm k_1 \dots  \pm k_l$ (possibly zero).
\end{proof}

\subsection{The odd spin algebra}\label{sec:oddspin}
In this section we are going to make heavy use of the following simple observation:
\begin{prop}
	If $A,B,C,D$ are linear operators and $B$ commutes with $A,C,D$ and $D$ commutes with $A,B,C$, then
	\begin{equation}
		[A\circ B, C \circ D] = B \circ D \circ [A,C]
	\end{equation}
\end{prop}
Let us consider the Lie subalgebra $\g$ of the arithmetic operators which is generated (in the sense of taking successive commutators) by the arithmetic operators $B^{\pm k} = \overline{\pm 2^{k-1}-\mp 2^{k-2}}$, for $2 \le k \le N $, and $B^{\pm 1} = \overline{\pm 1}$.
When it is necessary to distinguish this algebra for different integers $N$, we will denote it by $\g_N$.
\begin{prop}\label{prop:sl2k}
	We have that $B^k$ and $B^{-k}$ generate a subalgebra $\mathfrak{sl}_2^k$ isomorphic to $\mathfrak{sl}_2(\C)$.
	The elements $B^{-k}, [B^{-k}, B^{k}], B^k$ form a standard $\mathfrak{sl}_2$-triple. We have
	\begin{equation}
		\begin{split}
			&[B^{-k}, B^k] = \tfrac{1}{2}(p_k - p_{k-1}) \quad \text{ for } k>1,\\
			&[B^{-1}, B^1] = p_1
		\end{split}
	\end{equation}
\end{prop}
\begin{proof}
	We compute
	\begin{equation}
		\begin{split}
			&[B^{-k}, B^k] = \overline{-2^{k-1}}\circ \overline{+2^{k-1}}\circ \overline{+2^{k-2}}\circ \overline{-2^{k-2}} - \overline{+2^{k-1}}\circ \overline{-2^{k-1}}\circ \overline{-2^{k-2}}\circ \overline{+2^{k-2}} \\
			&= \tfrac{1}{4}(\id + p_k)(\id - p_{k-1}) -\tfrac{1}{4}(\id - p_k)(\id + p_{k-1}) = \tfrac{1}{2}(p_k - p_{k-1})
		\end{split}
	\end{equation}
	and also
	\begin{equation}
		\begin{split}
			[B^{\mp k}, \tfrac{1}{2}(p_k - p_{k-1})] = \tfrac{1}{2}(\overline{\pm 2^{k-2}})\circ [\overline{\mp 2^{k-1}}, p_k] -  \tfrac{1}{2}(\overline{\mp 2^{k-1}})\circ [\overline{\mp 2^{k-2}}, p_{k-1}] = \mp 2\cdot B^{\mp k}
		\end{split}
	\end{equation}
\end{proof}

\begin{theorem}\label{thm:oddspingen}
	The Lie algebra $\g$ is isomorphic to $\mathfrak{so}(2N+1,\C)$, and the representation on $\S_{2N+1}$ is equivalent to the fundamental spin representation.
\end{theorem}
\begin{proof}
	The Lie algebra $\g$ is semisimple, since it is generated by the $\mathfrak{sl}_2^k$-subalgebras from Prop. \ref{prop:sl2k}.
	It also comes equipped with a $\Z$--grading, inherited from the one given in Prop. \ref{prop:zgrading}.
	The elements of gradation $0$ form a subalgebra $\g_0$.
	This is abelian, and thus forms a Cartan subalgebra of $\g$. 
	We have that $\g_0$ is spanned by $\tfrac{1}{2}(p_N - p_{N-1}), \dots, \tfrac{1}{2}(p_{2} - p_{1}), p_1$.
	The positively graded subalgebra is generated by taking commutators of $B^N, \dots, B^{2}, B^1$; these form simple root vectors.
	Thus the Dynkin diagram of $\g$ can be obtained by computing such commutators.
	We get
	\begin{equation}
		\begin{split}
			&[B^2, B^1] = \overline{+2} \circ [B^{-1}, B^1] = p_1 \circ \overline{+2}\\
			&[p_1  \circ \overline{+2}, B^1] = \overline{+2} \circ [p_1, B^1] = -2 \overline{+2} \circ \overline{+1} = -2\cdot \overline{+3}\\
			&[\overline{+3}, B^1] = 0
		\end{split}
	\end{equation}
	and
	\begin{equation}
		\begin{split}
			&[B^k, B^{k-1}] = \overline{2^{k-1}}\circ\overline{2^{k-3}}\circ [\overline{-2^{k-2}},\overline{+2^{k-2}}] = p_{k-1} \circ \overline{+2^{k-1}-2^{k-3}}\\
			&[B^{k-1}, p_{k-1} \circ \overline{+2^{k-1}-2^{k-3}}] = (\overline{-2^{k-3}})^2 \circ (\dots) = 0\\
			&[B^{k}, p_{k-1} \circ \overline{+2^{k-1}-2^{k-3}}] = (\overline{-2^{k-1}})^2 \circ (\dots) = 0.
		\end{split}
	\end{equation}
	Here we are only interested in how many brackets between simple root vectors are nonzero, as this encodes the Dynkin diagram.
	We see that the Dynkin diagram is connected, so $\g$ is simple, and there is one short simple root corresponding to $B^1$, if $N>1$, since $B^2$ and $B^1$ admit two successive nonzero brackets, and all other simple roots are long.
	This is the Dynkin diagram of $\mathfrak{so}(2N+1,\C)$, and we have shown the first claim.\\

	For the second claim, we note that because of the $\Z$--grading, the element $[11\dots1]_2 \in \S_{2N+1}$ is the unique highest weight vector. 
	This vector has eigenvalue 1 for all parity operators, and therefore it has eigenvalue zero for $\tfrac{1}{2}(p_k - p_{k-1})$, but eigenvalue 1 for $p_1$.
	As $p_1$ corresponds to the short simple root vector $B^1$, and the corresponding root takes value 2 on $B^{-1}$, we see that $\S_{2N+1}$ has the same highest weight as the fundamental spin module, thus they are equivalent.
\end{proof}

\begin{cor}\label{cor:oddbasis}
	The Lie algebra $\g$ has basis elements 
	\begin{equation}
		\begin{split}
			&p_{l+1} \circ p_{l+2} \circ\dots \circ p_{k-1}\circ \overline{+2^k \pm 2^{l}},\\
			&p_i,\\
			&p_{l+1} \circ p_{l+2} \circ\dots \circ p_{k-1}\circ \overline{-2^k \pm 2^{l}},
		\end{split}
	\end{equation}
	where $0\le k<N$, $0<l<k$ and $0<i\le N$.
\end{cor}
\begin{proof}
	The proof will be by induction on $N$. First, if $N=1$, we have $\g_1=\langle \overline{+1}, p_1, \overline{-1} \rangle$, which establishes the base case.
	Suppose the result holds for $\g_{N-1}$. We will show that it holds for $\g_N$.
	The algebra is generated by the operators $\overline{+2^{k}-2^{k-1}}$, $\overline{-2^{k}+2^{k-1}}$. 
	The new generators compared to $\g_{N-1}$ are  $\overline{+2^{N-1}-2^{N-2}}$, $\overline{-2^{N-1}+2^{N-2}}$.
	Computing the commutators
	{\scriptsize
	\begin{equation}
		\begin{split}
			&[\overline{+2^{N-1}-2^{N-2}}, \overline{+2^{N-2}-2^{N-3}}] = \overline{+2^{N-1}}\circ \overline {-2^{N-3}}\circ [\overline{-2^{N-2}, +2^{N-2}}] = p_{N-1}\circ \overline{+2^{N-1}-2^{N-3}}\\
			&[p_{N-1}\circ \overline{+2^{N-1}-2^{N-3}}, \overline{+2^{N-3}-2^{N-4}}] = p_{N-1}\circ \overline{+2^{N-1}}\circ \overline {-2^{N-4}}\circ [\overline{-2^{N-3}, +2^{N-3}}] = p_{N-1}p_{N-2}\circ \overline{+2^{N-1}-2^{N-4}}\\
			&\vdots
		\end{split}
	\end{equation}
	} 
	leaves us with $p_{N-1}p_{N-3}\dots p_1 \circ \overline{+2^{N-1}}$ after $N-1$ steps.
	From here we compute the following commutators:
	\begin{equation}
		\begin{split}
			&[\overline{+1},p_{N-1}p_{N-2}\dots p_1 \circ \overline{+2^{N-1}}] = p_{N-1}\dots p_2 \circ \overline{+2^{N-1}}\circ[\overline{+1}, p_1] =  p_{N-1}\dots p_2 \circ \overline{+2^{N-1}+1}\\
			&[p_1\circ \overline{+2}, p_{N-1}\dots p_1 \circ \overline{+2^{N-1}}] = p_{N-1}\dots p_1^2 \circ \overline{+2^{N-1}}[\overline{+2}, p_2] = p_{N-1}\dots p_3 \circ \overline{+2^{N-1}+2}\\
			&\vdots\\
			&[p_1 p_2 \cdots p_{N-2}\circ \overline{+2^{N-2}}, p_{N-1}p_{N-3}\dots p_1 \circ \overline{+2^{N-1}}] = \overline{+2^{N-1} + 2^{N-2}}.
		\end{split}
	\end{equation}
	This has given us a total of $2N -1$ new positively graded elements (including the new generator). 
	The analogous computation gives $2N-1$ negatively graded elements.
	Finally we get $p_N$ from the commutator of $\overline{+2^{N-1} - 2^{N-2}}$ and $\overline{-2^{N-1} + 2^{N-2}}$.
	We know the difference in dimensions of $\g_N$ and $\g_{N-1}$ is $4N-1$, due to Theorem \ref{thm:oddspingen}, and a count shows that we have generated enough new elements. 
	These elements can be seen to all be linearly independent.
	Thus the formula holds for $\g_N$, and this concludes the proof.
\end{proof}

We can decompress the expressions from Cor. \ref{cor:oddbasis} into
{\scriptsize
\begin{equation}\label{eq:decompress1}
	\begin{split}
		&\overline{+1}\\
		&\overline{+2-1}\,, p_1\circ\overline{+2}\,, \overline{+2+1}\\
		&\overline{+4-2}\,, p_2\circ \overline{+4-1}\,, p_1 p_2\circ \overline{+4}\,, p_2\circ \overline{+4+1}\,,\overline{+4+2}\\
		&\vdots\\
		&\overline{+2^{N-1}-2^{N-2}}\,, p_{N-1} \circ \overline{+2^{N-1}-2^{N-3}}\,,\dots\,, p_1\cdots  p_{N-1} \circ \overline{+2^{N-1}}\,,\dots \,, p_{N-1} \circ \overline{+2^{N-1}+2^{N-3}}\,,\overline{+2^{N-1}+2^{N-2}}
	\end{split}
\end{equation}
}
and
{\scriptsize
\begin{equation}\label{eq:decompress2}
	\begin{split}
		&\overline{-1}\\
		&\overline{-2-1}\,, p_1\circ\overline{-2}\,, \overline{-2+1}\\
		&\overline{-4-2}\,, p_2\circ \overline{-4-1}\,, p_1 p_2\circ \overline{-4}\,, p_2\circ \overline{-4+1}\,,\overline{-4+2}\\
		&\vdots\\
		&\overline{-2^{N-1}-2^{N-2}}\,, p_{N-1} \circ \overline{-2^{N-1}-2^{N-3}}\,,\dots\,, p_1\cdots  p_{N-1} \circ \overline{-2^{N-1}}\,,\dots \,, p_{N-1} \circ \overline{-2^{N-1}+2^{N-3}}\,,\overline{-2^{N-1}+2^{N-2}}
	\end{split}
\end{equation}
}
together with the parity operators 
\begin{equation}\label{eq:decompress3}
	p_1, \dots, p_N.
\end{equation}

\subsection{The even spin algebra}
With our explicit description of the odd spin algebras from Cor. \ref{cor:oddbasis}, we can describe the even spin algebra, isomorphic to $\so(2N,\C)$, as its subalgebra.
Unfortunately it is not a graded subalgebra:
The even spin algebra $\d_N \simeq \so(2N,\C)$ is generated by adjoining the extra element $p_{N-1}p_{N-2}\dots p_1 \circ (\overline{+2^{N-1}} + \overline{-2^{N-1}})$ to $\g_{N-1}\subset \g_{N}$.
We caution that this is the only place so far where we have needed a linear combination of operators, as opposed to the formal linear combinations under the overlines.
\begin{theorem}\label{thm:evenspin}
	The even spin algebra $\d_{N}$ is the subalgebra of $\g_{N}$ generated by $\g_{N-1}$ together with the element $p_{N-1}p_{N-2}\dots p_1 \circ (\overline{+2^{N-1}} + \overline{-2^{N-1}})$
	It has a basis consisting of the elements coming from Cor. \ref{cor:oddbasis} applied to $\g_{N-1}$, in addition to the elements
	\begin{equation}
		\begin{split}
			\overline{\pm 2^{N-2}}\circ (\overline{+2^{N-1}} + \overline{-2^{N-1}})\\
			p_{N-1} \circ \overline{\pm 2^{N-3}}\circ (\overline{+2^{N-1}} + \overline{-2^{N-1}})\\
			p_{N-1}p_{N-2} \circ \overline{\pm 2^{N-4}}\circ (\overline{+2^{N-1}} + \overline{-2^{N-1}})\\
			\vdots\\
			p_{N-1}\cdots p_2 \circ \overline{\pm 1}\circ (\overline{+2^{N-1}} + \overline{-2^{N-1}})
		\end{split}
	\end{equation}
\end{theorem}
\begin{proof}
	Let us denote $(\overline{+2^{N-1}} + \overline{-2^{N-1}})=C$, and note that this linear operator commutes with $\g_{N-1}$.
	The indicated elements can be generated by taking commutators between $p_{N-2}p_{N-3}\dots p_1 \circ C$ and the elements
	\begin{equation}
		\begin{split}
			\overline{\pm 1},\, p_1 \circ\overline{\pm 2},\, p_1p_2 \circ\overline{\pm 4}, \dots, p_1p_2\cdots p_{N-2}\circ \overline{\pm 2^{N-2}}
		\end{split}
	\end{equation}
	from $\g_{N-1}$.
	This yields $2N-1$ new elements.
	The elements $\pm1, p_1,\, p_{N-1}\cdots p_1 C, p_{N-1}\cdots p_2 \circ \overline{\pm 1}  C$, form a subalgebra isomorphic to $\sl(2,\C)\oplus \sl(2,\C)$.
	Thus the whole Lie algebra can be seen to be generated by non-commuting copies of $\sl(2,\C)$. 
	Therefore the Lie algebra generated is semi-simple, and since it contains $\g_{N-1}$ and is contained in $\g_N$, it must be $\d_N$
	We also get that $p_{N-1}\dots p_1 \circ C$ together with $p_1, \dots, p_{N-2}$ forms a Cartan subalgebra.
\end{proof}

\section{Graph algorithm}
Another way to interpret the generators $B^{\pm k}$ from Section \ref{sec:oddspin} is as a ``quantum version'' of bit-shift operators, in light of the interpretation of $\S_{2N+1}$ as the $N$-qubit state space.
Then we may think of $B^k$ as shifting states with nonzero $k-1$st bit and zero $k$th bit to states with nonzero $k$th bit and zero $k-1$st bit.
For example,
\begin{equation}
	B^2([01]_2) = [10]_2
\end{equation}
while the operators $B^{\pm 1}$ creates a new nonzero bit at the left or right edge, if that is possible.
\begin{equation}
	B^1([00]_2) = [01]_2
\end{equation}
This interpretation makes it possible to encode the generator structure of $\so(2N+1,\C)$ and its spin representation in a coloured graph, by writing the action of all possible bitshift operators.
\begin{example}
	The positively graded generators of the spin matrices of $\mathfrak{spin}(7,\C)$ are encoded in the following diagram:
	\[\begin{tikzcd}
		{[000]} & {[001]} & {[010]} & {[011]} \\
		& {} & {[100]} & {[101]} & {[110]} & {[111]}
		\arrow["{B^1}", from=1-1, to=1-2]
		\arrow["{B^2}", from=1-2, to=1-3]
		\arrow["{B^1}", from=1-3, to=1-4]
		\arrow["{B^3}", from=1-3, to=2-3]
		\arrow["{B^3}", from=1-4, to=2-4]
		\arrow["{B^1}", from=2-3, to=2-4]
		\arrow["{B^2}", from=2-4, to=2-5]
		\arrow["{B^1}", from=2-5, to=2-6]
	\end{tikzcd}\]
	The nonzero matrix coefficients can be read off by considering binary numbers $[m]$ as basis vectors in $\S_{2N+1}$, and the  negatively graded generators are their transposes.
	The whole spin representation is obtained by taking commutators.
\end{example}
We offer one further example:
\begin{example}
	We consider the case the spin representation of $\so(9,\C)$. Here we compactify the diagram by encoding the different generators by colours (pink,blue,green,orange), and consider each node a distinct basis vector in a free complex vector space.
	
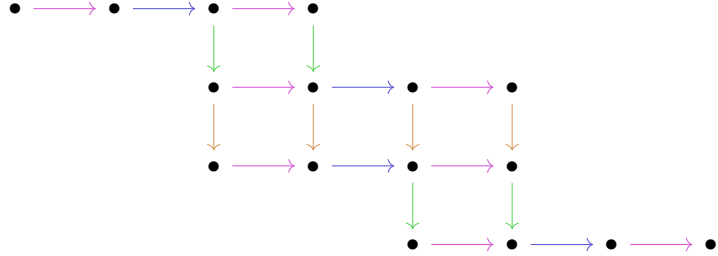
\[\begin{tikzcd}
		\bullet & \bullet & \bullet & \bullet \\
		&& \bullet & \bullet & \bullet & \bullet \\
		&& \bullet & \bullet & \bullet & \bullet \\
		&&&& \bullet & \bullet & \bullet & \bullet
		\arrow[color={rgb,255:red,214;green,92;blue,214}, from=1-1, to=1-2]
		\arrow[color={rgb,255:red,92;green,92;blue,214}, from=1-2, to=1-3]
		\arrow[color={rgb,255:red,214;green,92;blue,214}, from=1-3, to=1-4]
		\arrow[color={rgb,255:red,92;green,214;blue,92}, from=1-3, to=2-3]
		\arrow[color={rgb,255:red,214;green,153;blue,92}, from=2-3, to=3-3]
		\arrow[color={rgb,255:red,92;green,214;blue,92}, from=1-4, to=2-4]
		\arrow[color={rgb,255:red,214;green,153;blue,92}, from=2-4, to=3-4]
		\arrow[color={rgb,255:red,214;green,92;blue,214}, from=3-3, to=3-4]
		\arrow[color={rgb,255:red,92;green,92;blue,214}, from=3-4, to=3-5]
		\arrow[color={rgb,255:red,92;green,92;blue,214}, from=2-4, to=2-5]
		\arrow[color={rgb,255:red,214;green,92;blue,214}, from=2-3, to=2-4]
		\arrow[color={rgb,255:red,214;green,92;blue,214}, from=2-5, to=2-6]
		\arrow[color={rgb,255:red,214;green,92;blue,214}, from=3-5, to=3-6]
		\arrow[color={rgb,255:red,214;green,153;blue,92}, from=2-6, to=3-6]
		\arrow[color={rgb,255:red,214;green,153;blue,92}, from=2-5, to=3-5]
		\arrow[color={rgb,255:red,92;green,214;blue,92}, from=3-5, to=4-5]
		\arrow[color={rgb,255:red,92;green,214;blue,92}, from=3-6, to=4-6]
		\arrow[color={rgb,255:red,214;green,92;blue,214}, from=4-5, to=4-6]
		\arrow[color={rgb,255:red,92;green,92;blue,214}, from=4-6, to=4-7]
		\arrow[color={rgb,255:red,214;green,92;blue,214}, from=4-7, to=4-8]
	\end{tikzcd}\]	
	The arrows encode nonzero matrix coefficients as above, and $\mathfrak{spin}(9,\C)$ is obtained by taking these matrices, their transposes, and a sufficient number of commutators.
\end{example}
We note that it is easily possible to see lower dimensional spin algebras as subdiagrams, and hence infer branching rules, for example.


\section{Compact real form $\Spin(N)$}\label{sec:compactspin}
We are often interested in the compact real form $\so(2N+1)$ rather than $\so(2N+1,\C)$.
Since we have the grading structure of $\g_N$, going to the compact real form can be done as follows:
\begin{prop}
	The real Lie algebra generated by the operators
	\begin{equation}
		\begin{split}
			&Q^k_+ = B^k-B^{-k}\\
			&Q^k_- = i\cdot(B^k+B^{-k})\\
			&p_l^q = i\cdot p_l
		\end{split}
	\end{equation}
	on $\S_{2N+1}$, where $i$ is the imaginary unit, is isomorphic to $\so(2N+1)$, and the representation is the fundamental complex spin representation.
\end{prop}
\begin{proof}
	Follows from general structure theory.
\end{proof}
One can get all the other basis elements of the compact form by taking all pairs of basis elements $P\circ \overline{\pm 2^k \mp 2^l}$ from $\g_N$, which are related by a sign change and where $P$ is some product of parity operators, and taking the combinations
\begin{equation}
	\begin{split}
		P \circ (\overline{\pm 2^k \mp 2^l} - \overline{\mp 2^k \pm 2^l})\\
		i\cdot P \circ (\overline{\pm 2^k \mp 2^l} + \overline{\mp 2^k \pm 2^l})
	\end{split}
\end{equation}
and similarly for $P\circ \overline{\pm 2^k \pm 2^l}$.

\subsection{Real structure}
However, the spin representation also sometimes admits a real structure, i.e. a basis where all representation matrices are real.
When this exists, it must coincide with the basis in which all the matrices of the Cartan subalgebra are real.
We can describe this basis in our terms.
\begin{prop}
	Consider the involution $\tau$ on the set of binary numbers of fixed length $N$ given by flipping all bits.
	Then the operators $p_l^q$ have real matrix coefficients in the basis
	\begin{equation}
		([0]+i\cdot\tau([0]),\, i\cdot [0]+\tau([0]),\, [1]+i\cdot\tau([1]),\, i\cdot [1]+\tau([1]),\, \dots,\, i\cdot [2^{N-1}-1]+\tau([2^{N-1}-1]) ).
	\end{equation}
\end{prop}
\begin{proof}
	The involution $\tau$ takes weight vectors with weight $\lambda$ to weight vectors with weight $-\lambda$. 
	Thus each $p_l^q$ will act as a rotation on the real plane spanned by $[m]+i\cdot\tau([m])$ and $i\cdot [m]+\tau([m])$.
\end{proof}

\section{Future directions}
It seems that the generators and their relations do not really depend strongly on $N$. Thus one could extend $\{1,0\}^N$ to $\{1,0\}^\mathbb{N}$ and introduce an infinite dimensional ``spinfinity'' algebra with operators given by the same formulas as in the finite dimensional case.

\end{document}